\documentclass{article}
\usepackage{amsrefs}
\usepackage{amssymb}
\usepackage{amsmath}
\usepackage{amsthm}
\usepackage{graphicx}
\usepackage{ytableau}

\newtheorem{theorem}{Theorem}% [chapter]
\newtheorem{proposition}{Proposition}%[chapter]
\newtheorem{lemma}{Lemma}%[chapter]
\newtheorem*{corollary}{Corollary}
\newtheorem*{examples}{Theorem B}
\newtheorem*{counterexamples}{Theorem C}
\newtheorem*{Martens}{Theorem A (H. Martens' Theorem for chains of cycles)}

\theoremstyle{definition}
\newtheorem{definition}{Definition}%[chapter]

\theoremstyle{remark}
\newtheorem{example}{Example}%[chapter]

\DeclareMathOperator{\im}{im}
\DeclareMathOperator{\divi}{div}
\DeclareMathOperator{\rank}{rk}

\begin{document}

\author{Marc Coppens\footnote{KU  Leuven; Department of Mathematics, Section of Algebra,
Celestijnenlaan 200B bus 2400 B-3001 Leuven, Belgium; Department of Elektrotechniek (ESAT), Technologiecampus Geel, Kleinhoefstraat 4, 2440 Geel, Belgium; email: marc.coppens@kuleuven.be.}}

\title{A study of H. Martens' Theorem on chains of cycles}
\date{}

\newtheorem*{subject}{2000 Mathematics Subject Classification}
\newtheorem*{keywords}{Keywords}

\maketitle \noindent

\begin{abstract}
Let $C$ be a smooth curve of genus $g$ and let $d$,$r$ be integers with $1 \leq r \leq g-2$ and $2r\leq d \leq g-2+r$.
H. Martens' Theorem states that $\dim (W^r_d(C))=d-2r$ implies $C$ is hyperelliptic.
It is known that for a metric graph $\Gamma$ of genus $g$ such statement using $\dim (W^r_d(\Gamma))$ does not hold.
However replacing $\dim (W^r_d(\Gamma))$ by the so-called Brill-Noether rank $w^r_d(\Gamma)$ it was stated as a conjecture.
Using a similar definition in the case of curves one has $\dim (W^r_d(C))=w^r_d(C)$.

Let $\Gamma$ be a chain of cycles of genus $g$ and let $r,d$ be integers with $1 \leq r \leq g-2$ and $2r\leq d \leq g-3+r$.
If $w^r_d(\Gamma)=d-2r$ then we prove $\Gamma$ is hyperelliptic.
In case $g \geq 2r+3$ then we prove there exist non-hyperelliptic chains of cycles satisfying $w^r_{g-2+r}(\Gamma)=g-2-r$, contradicting the conjecture.
We give a complete description of all counterexamples within the set of chains of cycles to the statement of H. Martens' Theorem.
Those counterexamples also give rise to chains of cycles such that $w^r_{g-2+r}(\Gamma) \neq w^1_{g-r}(\Gamma)$.
This shows that the Riemann-Roch duality does not hold for the Brill-Noether ranks of metric graphs.
\end{abstract}

\begin{subject}
05C25, 14T15
\end{subject}

\begin{keywords}
Brill-Noether ranks, divisors, hyperelliptic graphs, metric graphs
\end{keywords}

\section{Introduction}\label{section1}

During the recent decades, a theory of divisors on metric graphs is developed having lots of properties similar to those of the theory of divisors on smooth projective curves. For the motivations coming from 1-parameter degenerations of families of smooth projective curves and how it occurs as tropicalization of divisors in Berkovich analytic theory, we refer to \cite{ref2}.

In particular,  if $D$ is a divisor on a metric graph $\Gamma$ (see Section \ref{section2}) one defines the rank $\rank (D)$ (see Definition \ref{def3}).
As a topological space $\Gamma$ has a genus $g(\Gamma)$ and there is a Riemann-Roch Theorem as in the case of smooth curves (see \cite{ref3}, \cite{refMZ}).
As in the case of smooth curves this Riemann-Roch Theorem gives a lower bound $\rank (D) \geq \deg(D)-g(\Gamma)$.
In case $\deg (D)>2g(\Gamma)-2$ this always is an equality, so one is interested in $\rank (D)$ for $\deg (D) \leq 2g(\Gamma)-2$.

If $C$ is a smooth curve of genus $g$ one has Clifford's Theorem giving an upper bound on $r(D)=\dim \vert D \vert$ and a description of the divisors attaining that upper bound. In particular it implies $r(D) \leq \deg (D)/2$ in case $\deg (D) \leq 2g(C)-1$ and in case $D$ is not equivalent to $0$ or $K_C$ then $r(D)=\deg (D)/2$ implies $C$ is a hyperelliptic curve. The inequality $r(D) \leq \deg (D)/2$ in the case of metric graphs is proved in \cite{ref4} while in \cite{ref5} it is proved that $r(D)= \deg (D)/2$ in case $D$ is not equivalent to $0$ or $K_{\Gamma}$ implies $\Gamma$ is a hyperelliptic graph.

In the case of curves there is an important generalization of Clifford's Theorem called H. Martens' Theorem.
It says that if $C$ is a smooth curve of genus $g$, $r$ is an integer with $0 \leq r \leq g-1$ and $d$ is an integer with $2r \leq d \leq g-1+r$ then $\dim (W^r_d(C)) \leq d-2r$.
Moreover in case $1 \leq r \leq g-2$ and $2r \leq d \leq g-2+r$ then $W^r_d(C)$ has a component of dimension $d-2r$ if and only if $C$ is hyperelliptic (see \cite{ref7}).
Of course the inequality $d \geq 2r$ comes from Clifford's Theorem.
It would be nice if a similar statement does hold for graphs.

It should be noted that because of the Riemann-Roch Theorem we have $\dim (W^r_d(C))=\dim (W^{r-d+g-1}_{2g-2-d}(C))$ and therefore it is enough to consider the case $d \leq g-1$.
In particular, in \cite{ref6} H. Martens's Theorem is stated under the condition $2 \leq d \leq g-1$ and $r$ is an integer with $0 < 2r \leq d$.

In case $\Gamma$ is a metric graph then as a set $W^r_d(\Gamma)$ is equal to the set of equivalence classes of divisors $D$ on $\Gamma$ of degree $d$ and rank at least $r$.
In \cite{ref8} the authors define a polyhedral set structure on $W^r_d(\Gamma)$ such that we can speak of  $\dim (W^r_d(\Gamma))$.
However in the same paper the authors consider an example of a metric graph $\Gamma$ satisfying $\dim (W^1_3(\Gamma))=1$ while $\Gamma$ is not a hyperelliptic graph. This contradicts the H. Martens' Theorem using $\dim (W^r_d (\Gamma))$.
Also the behaviour in families of metric graphs $\Gamma_t$ of $\dim (W^r_d(\Gamma_t))$ is not nice: it does not vary upper-semicontinuously.

In that paper \cite{ref8} the authors introduce a Brill-Noether rank $w^r_d(\Gamma)$ (see Definition \ref{def3}) such that for a similar definition using curves $C$ one has $w^r_d(C)=\dim(W^r_d(C))$.
In particular from the Riemann-Roch Theorem it follows that $w^r_d(C)=w^{r-d+g-1}_{2g-2-d}(C)$.
So, in the proof of H. Martens' Theorem for curves one can restrict to $d \leq g-1$.
In \cite{ref8} the authors show that in case of their above-mentioned example of a non-hyperelliptic graph with $\dim (W^1_3(\Gamma))=1$, one has $w^1_3(\Gamma)=0$.
So that example does not contradict H. Martens' Theorem for graphs using the Brill-Noether ranks $w^r_d$.
Also it is proved that those Brill-Noether ranks vary upper-semicontinuously in families of graphs.
So it is natural to consider H. Martens' Theorem for metric graphs using those Brill-Noether ranks $w^r_d(\Gamma)$.
As such it is stated as a conjecture in the appendix of \cite{ref10} written by Y. Len and D. Jensen but only under the conditions $d \leq g-1$ and $0< 2r \leq d$.
In that appendix they prove the inequality $w^r_d(\Gamma)\leq d-2r$ (under the assumption of the above mentioned inequalities on $d$ and $r$).
Unlike the case of curves it is not a priori clear that the Riemann-Roch Theorem for graphs would imply $w^r_d(\Gamma)=w^{r-d+g-1}_{2g-2-d}(\Gamma)$.
So there is no a priori reason why to restrict the statement of H. Martens' Theorem for graphs to the cases $2 \leq d \leq g-1$ with $0 < 2r \leq d$.
So in this paper we will consider H. Martens' Theorem for graphs under the assumptions $1 \leq r \leq g-2$ and $2r \leq d \leq g-2+r$.
(Using the arguments from \cite{ref10} one finds that the inequality $w^r_d(\Gamma) \leq d-2r$ holds for $0 \leq r \leq g-1$ with $2r \leq d \leq g-1+r$.)

A chain of cycles is a special type of metric graph (see Section \ref{section2}) that entered in the scene in a paper giving a tropical proof of classical Brill-Noether Theory (see \cite{ref9}).
Afterwards those graphs where used to obtain a lot of results concerning smooth curves (see e.g. \cite{ref11}).
In \cite{ref1} a very handsome description of the spaces $W^r_d(\Gamma)$ in the case of chains of cycles is given (see Theorem \ref{theorem1}).
This was used in \cite{ref11} amongst others and it makes it possible to use chain of cycles to get a better idea of what can happen in the theory of divisors on metric graphs.
As an example, in \cite{ref12} one studies the scrollar invariants of equivalence classes of divisors of rank 1 and one finds similarities but also differences with the theory of scrollar invariants on curves.
In this paper we study H. Martens' Theorem in the case of chains of cycles.

We obtain the following theorem.
\begin{Martens}\label{theoremA}
Let $\Gamma$ be a chain of cycles of genus $g$.
Let $r$ be an integer with $1 \leq r \leq g-2$ and let $d$ be an integer with $2r \leq d \leq g-3+r$.
Then $w^r_d(\Gamma)=d-2r$ implies $\Gamma$ is hyperelliptic.
\end{Martens}
Compared to the statement of H. Martens' Theorem for curves, this Theorem A misses the case $d=g-2+r$.
Compared to the conjecture stated in \cite{ref10} (using the inequality $d\leq g-1$ in its statement) it misses the case $(d,r)=(g-1,1)$.
As a matter of fact, the methods of the proof of Theorem A only shows that the cases $w^r_{g-2+r}(\Gamma)=g-2-r$ and $r+2 \leq g \leq 2r+2$ imply $\Gamma$ is hyperelliptic.
In particular for $d = g-1$ and $w^1_{g-1}(\Gamma)=g-3$ only $g\in \{3,4\}$  implies $\Gamma$ is hyperelliptic.
In general the methods of the proof only imply that $\Gamma$ is a very particular kind of chain of cycles in case $w^r_{g-2+r} (\Gamma)=g-2-r$.
Those non-hyperelliptic chains of cycles will be called Martens-special chains of cycles for rank $r$ (and they also have a certain type $k$; see Definition \ref{def5}).
From this definition it follows that in case $r \geq 2$ a Martens-special chain of cycles for rank $r$ also is a Martens-special chain of cycles for rank $r-1$.
In particular , a Martens-special chain of cycles for rank 1 is simply called a Martens-special chain of cycles.
\begin{examples}\label{theoremB}
A chain of cycles $\Gamma$ satisfies $w^r_{g-2+r}(\Gamma)=g-2-r$ for some integer $1 \leq r \leq g-2$ if and only if either $\Gamma$ is hyperelliptic or $\Gamma$ is a Martens-special chain of cycles for rank $r$.
\end{examples}
In Lemma \ref{lemma6} we prove that a Martens-special chain of cycles of type $k$ has gonality at most $k+2$.
This implies a weak H. Martens' Theorem for $w^r_{g-2+r}(\Gamma)$ for chains of cycles (see Proposition \ref{prop2}).

Since there exist Martens-special chains of cycles for rank $r$ for each genus $g\geq 2r+3$ this shows that the statement of H. Martens' Theorem for curves does not hold in its full generality for metric graphs.
Combining Theorems A and B we also obtain the following theorem. 

\begin{counterexamples}
The equality between $w^r_d(\Gamma)$ and $w^{r-d+g-1}_{2g-2-d}(\Gamma)$ does not hold in its full generality for metric graphs.
\end{counterexamples}
This confirms that there is no reason to restrict to the case $d\leq g-1$ in considering H. Martens' Theorem for metric graphs.

In Section \ref{section2} we recall some generalities and we give a summary of some crucial results from \cite{ref1} used in this paper.
In Section \ref{section3} we prove Theorem A and we obtain the definition of Martens-special chains of cycles for rank $r$.
Those Martens-special chains of cycles are studied in more detail in Section \ref{section4}, proving Theorems B and C.

\section{Generalities}\label{section2}

A metric graph $\Gamma$(shortly a graph) is a compact connected metric space such that for each point $P$ on $\Gamma$ there exists $n\in \mathbb{Z}$ with $n\geq 1$ and $r\in \mathbb{R}$ with $r>0$ such that some neighbourhood of $P$ is isometric to $\{ z \in \mathbb{C} : z=te^{2\pi ik/n} \text{ with } t\in [0,r]\subset \mathbb{R} \text{ and } k\in \mathbb{Z} \}$ and with $P$ corresponding to $0$.
In case $n\neq 2$ then we say $P$ is a vertex of $\Gamma$.
After omitting the vertices from $\Gamma$, the resulting connected components are called the edges of $\Gamma$.

\begin{figure}[h]
\begin{center}
\includegraphics[height=2 cm]{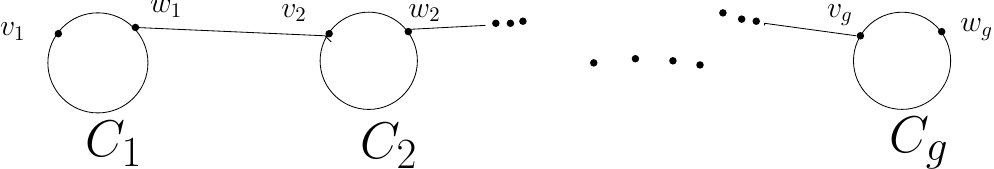}
\caption{a chain of cycles }\label{Figuur 1}
\end{center}
\end{figure}

In this paper the graph $\Gamma$ is a chain of cycles.
Such a chain of cycles is obtained as follows.
We take $g$ graphs $C_1, \cdots , C_g$ ($g \in \mathbb{Z}$ with $g \geq 1$) each of them isometric to a circle (such graphs are called cycles).
On each cycle we choose two different points $v_i$ and $w_i$ and we take an orientation on the cycle such that the positively oriented arc segment from $v_i$ to $w_i$ is the shortest one (so we fix one of the two possibilities when both arcs have equal length).
For $1 \leq i \leq g-1$ we connect $w_i$ to $v_{i+1}$ using a graph isometric to a straight line segment.
The resulting graph is called a chain of cycles of genus $g$ (see Figure \ref{Figuur 1}).
We will use those notations ($C_i, v_i, w_i$) for a chain of cycles throughout the paper.

In principle the length of the straight line segments can be choosen to be zero, but for distinguishing points on $C_i$ and $C_{i+1}$ in case $1 \leq i \leq g-1$ it is handsome to have $w_i \neq v_{i+1}$.
Moreover such bridges don't affect the Brill-Noether theory on the graph, so it is harmless to contract or de-contract them.
So we assume those lengths are non-zero.

Such a chain of cycles $\Gamma$ has a torsion profile $\underline{m}=(m_2, \cdots , m_{g-1})$ consisting of $g-2$ non-negative integers whose definition comes from \cite{ref1}, Definition 1.9.

\begin{definition}\label{def1}
Let $\ell_i$ be the length of the cycle $C_i$ and let $\ell(v_i,w_i)$ be the length of the positively oriented arc on $C_i$ from $v_i$ to $w_i$.
In case $\ell_i$ is an irrational multiple of $\ell(v_i,w_i)$ then $m_i=0$.
Otherwise $m_i$ is the minimum positive integer such that $m_i \cdot \ell(v_i,w_i)$ is an integer multiple of $\ell_i$.
\end{definition}

Following \cite{ref1}, Definition 3.1, we use the following convention to denote points on $C_i$.

\begin{definition}\label{def2}
For $\xi \in \mathbb{R}$ let $\langle\xi \rangle_i$ denote the point on $C_i$ that is located $\xi \cdot\ell(v_i,w_i)$ units on the positive oriented path on $C_i$ starting at $w_i$.
\end{definition}

It follows that, for integers $n_1$,$n_2$ one has $\langle n_1\rangle_i=\langle n_2\rangle_i$ if and only if $n_1 \equiv n_2 \mod{m_i}$ ($m_i$ as in Definition \ref{def1}).

A divisor $D$ on an arbitrary graph $\Gamma$ is a finite formal linear combination $\sum_{P \in \Gamma}n_P\cdot P$ of points on $\Gamma$ with integer coefficients (hence $n_P \neq 0$ for only finitely many points).
The degree $\deg (D)$ of $D$ is $\sum _{P \in \Gamma}n_P$.
We say that $D$ is effective if $n_P\geq 0$ for all $P \in \Gamma$.
A rational function on $\Gamma$ is a continuous function $f:\Gamma \rightarrow \mathbb{R}$ that can be described as a piecewise affine function with integer slopes on the edges.
For $P\in \Gamma$ we define $f_P$ as being the sum of all slopes of $f$ on $\Gamma$ in all directions emanating from $P$.
In this way $f$ defines a divisor $\divi (f)=\sum _{P\in \Gamma}f_P\cdot P$.

\begin{definition}\label{def3}
Two divisors $D_1$ and $D_2$ are called equivalent if $D_2-D_1=\divi (f)$ for some rational function $f$ on $\Gamma$.

For a divisor $D$ we define the rank $\rank (D)$ as follows.
In case $D$ is not equivalent to an effective divisor then $\rank (D)=-1$.
Otherwise, $\rank (D)$ is the maximal integer $r$ such that for each effective divisor $E$ of degree $r$ on $\Gamma$ there exists an effective divisor $D'$ on $\Gamma$ equivalent to $D$ and containing $E$ (meaning $D'-E$ is effective).

The Brill-Noether number $w^r_d(\Gamma)$ is -1 in case $\Gamma$ contains no divisor of degree $d$ and rank at least $r$.
Otherwise it is the largest integer $w$ such that for every effective divisor $E$ of $\Gamma$ of degree w+r there exists an effective divisor $D$ of degree $d$ and rank at least $r$ containing $E$.
\end{definition}

In case $\Gamma$ is a chain of cycles of genus $g$ then any rational function on $C_i$ can be extended to a rational function on $\Gamma$ using constant functions on the connected components of $(\Gamma \setminus C_i)\cup \{ v_i, w_i \}$.
In particular if two divisors on $C_i$ are equivalent as divisors on $C_i$ then they also are equivalent as divisors on $\Gamma$.
In particular for any divisor $D$ of degree $d$ on $C_i$ there exist unique points $P'$ and $P''$ on $C_i$ such that on $\Gamma$ the divisor $D$ is equivalent to the divisors $(d-1)w_i+P'$ and $(d-1)v_i+P''$.
This follows from the Riemann-Roch Theorem on $C_i$.
It follows that every effective divisor of degree $d \leq g$ on $\Gamma$ is equivalent to an effective divisor consisting of $d$ points on $d$ different cycles $C_i$.
In particular in \cite{ref1}, Lemma 3.3 one obtains the following lemma.

\begin{lemma}\label{lemma1}
Let $D$ be any divisor of degree $d$ on a chain of cycles of genus $g$.
Then $D$ is equivalent to a unique divisor of the form $\sum_{i=1}^g \langle\xi _i\rangle_i +(d-g)\cdot w_g$ (here $\xi _i \in \mathbb{R}$ and $\langle\xi _i\rangle_i \in C_i$ using Definition \ref{def2}).
\end{lemma}

We are going to call this unique divisor from Lemma \ref{lemma1} the representing divisor of (the equivalence class of) $D$.

In \cite{ref1}, using those representing divisors, one gives a nice description for the set $W^r_d( \Gamma )$ of equivalence classes of divisors of degree $d$ and rank at least $r$.
We first use the following definition from \cite{ref1}.

\begin{definition}\label{def4}
For positive integers $m$ and $n$ we write $[m \times n]$ to denote the set $\{ 1, \cdots , m\}\times \{1, \cdots, n\}$.
It is represented by a rectangle with $n$ rows and $m$ columns (see e.g. Figure \ref{Figuur 2}).
The element on the $i$-th column and the $j$-th row corresponds to $(i,j) \in [m\times n]$ (the numbering is such that $(1,1)$ corresponds to the lower left corner in Figure \ref{Figuur 2}).

\begin{figure}
\begin{center}
\ydiagram{12 , 12, 12}
\caption{rectangle $[12 \times 3]$}\label{Figuur 2}

\end{center}
\end{figure}

An $\underline{m}$-displacement tableau on $[m \times n]$ of genus $g$ is a function $t:[m \times n] \rightarrow \{ 1, \cdots ,g\}$ such that
\begin{enumerate}
\item $t$ is strictly increasing in any given row and column of $[m \times n]$.
\item if $t(x,y)=t(x',y')$ then $x-y \equiv x'-y' \mod{m_{t(x,y)}}$.
\end{enumerate}
\end{definition}

\begin{example}\label{example1}
In case $\underline{m}=(2,3,0,3,2,2,2,3,5,2,3)$ Figure \ref{Figuur 3} shows an $\underline{m}$-displacement tableau $t$ on $[5 \times 3]$  by denoting the values in the corresponding boxes of the rectangle.
\begin{figure}[h]
\begin{center}
\begin{ytableau}
3 & 7 & 8 & 11 & 12\\
2 & 5 & 7 & 10 & 11\\
1 & 2 & 4 & 5 & 8
\end{ytableau}
\caption{Example \ref{example1} }\label{Figuur 3}
\end{center}
\end{figure}
\end{example} 

The following statement is part of Theorem 1.4 in \cite{ref1}.

\begin{theorem}\label{theorem1}
A divisor $D$ on $\Gamma$ of degree $d$ has rank at least $r$ if and only if on $\lambda = [(g-d+r) \times (r+1)]$ there exists an $\underline{m}$-displacement tableau $t$ such that for the representing divisor $\sum_{i=1}^g \langle\xi _i\rangle_i +(d-g)w_g$ of $D$ one has $\xi _{t(x,y)} \equiv x-y \mod{m_{t(x,y)}}$ for all $(x,y) \in \lambda$ (because of condition 2 in Definition \ref{def4} such divisor is well-defined).
\end{theorem}

Given an $\underline{m}$-displacement tableau $t$ on $\lambda =[(g-d+r) \times (r+1)]$ we write $\mathbb{T}(t)$ to denote the space of representing divisors $\sum_{i=1}^g \langle\xi _i\rangle_i +(d-g)w_g$ with $\xi _{t(x,y)} \equiv x-y \mod{m_{t(x,y)}}$ for all $(x,y)\in \lambda$.
Clearly $\mathbb{T}(t)$ has the structure of a torus of some dimension $m$ with $m$ the number of values in $\{ 1, \cdots ,g \} \setminus \im (t)$.
It can be identified with the product of those $m$ cycles $C_j$ with $j \notin \im (t)$.
In Example \ref{example1} we have $\mathbb{T}(t) \subset W^2_9(\Gamma)$ and $\dim ( \mathbb{T}(t))=2$ since $\mathbb{T}(t)$ can be identified with $C_6 \times C_9$.

The theorem implies $W^r_d(\Gamma)$ can be identified with the union of those spaces $\mathbb{T}(t)$ with $t$ ranging over the set of $\underline{m}$-displacement tableaux on $[(g-d+r) \times (r+1)]$.

\begin{example}\label{example2}
A very important example for this paper is the case when $\Gamma$ is hyperelliptic, meaning $\Gamma$ has a divisor of degree 2 and rank 1.
This is possible if and only if there is an $\underline{m}$-tableau of genus $g$ on $[(g-1) \times 2]$.
Because of the first condition of Definition \ref{def4} the only possibility is $t(x,1)=x$ and $t(x,2)=x+1$ for $1 \leq x \leq g-1$.
Because of the second condition of Definition \ref{def4} this is possible if and only if $m_i=2$ for $2 \leq i \leq g-1$.
This tableau is called the hyperelliptic tableau of genus $g$ (see Figure \ref{Figuur 4}).
\begin{figure}[h]
\begin{center}
\begin{ytableau}
2 & 3 & 4 & 5 & 6\\
1 & 2 & 3 & 4 & 5
\end{ytableau}
\caption{The hyperelliptic tableau of genus 6 }\label{Figuur 4}
\end{center}
\end{figure}
\end{example}

\section{H. Martens' Theorem for chains of cycles}\label{section3}

We start by given a very easy proof of the main result of \cite{ref5} in the case of chains of cycles.
It will be used to obtain the main ingredient in the proof of Theorem A.

\begin{proposition}\label{prop1}
Assume a chain of cycles $\Gamma$ of genus $g$ has a divisor $D$ of rank $r$ between 1 and $g-2$ and degree $2r$.
Then $\Gamma$ is a hyperelliptic chain of cycles.
\end{proposition} 
 
\begin{proof}
The representing divisor of $D$ corresponds to an $\underline{m}$-displacement tableau $t$ on $\lambda =[(g-r) \times (r+1)]$.
Take $(x,y) \in \lambda$ and consider the sequence of integers $t(1,1)<t(2,1)< \cdots t(x,1) <t(x,2) < \cdots <t(x,r+1) <t(x+1,r+1) < \cdots < t(g-r,r+1)$.
This sequence consists of $g$ different integers between $1$ and $g$, hence it is the sequence $1 < 2 < \cdots <g$.
It follows that $t(x,y)=x+y-1$.
If $(x,y) \neq (1,1)$ and $(x,y) \neq (g-r,r+1)$ then $(x-1,y+1)\in \lambda$ or $(x+1,y-1)\in \lambda$.
In those cases $t(x-1,y+1)$ or $t(x+1,y-1)$ is equal to $t(x,y)$.
This implies $m_2 = m_3 = \cdots = m_{g-1}=2$, hence $\Gamma$ is hyperelliptic.
\end{proof}

The proof of Proposition \ref{prop1} generalizes to show that a chain of cycles has a lot of indices $i$ with $m_i=2$ in case $\dim (W^r_d)=d-2r$ for some $1\leq r \leq g-2$.

\begin{lemma}\label{lemma2}
If $\Gamma$ is a chain of cycles of genus $g$ and $t$ is an $\underline{m}$-displacement tableau $t$ such that $\mathbb{T}(t) \subset W^r_d(\Gamma)$ then $\dim (\mathbb{T}(t)) \leq d-2r$.
In case $\dim (\mathbb{T}(t))=d-2r$ and $\im (t)=\{ i_1, \cdots ,i_{g-d+2r} \}$ with $i_1 < \cdots < i_{g-d+2r}$ then $m_{i_2}= \cdots = m_{i_{g-d+2r-1}}=2$.
In words, $\im(t)$ is a subset of $\{ 1, \cdots , g \}$ such that each $i \in \im (t)$ different from the smallest and the largest element of $\im (t)$ satisfies $m_i=2$.
\end{lemma}

\begin{proof}
Let $t$ be an $\underline{m}$-displacement tableau on $\lambda=[(g-d+r) \times (r+1)]$.
Choose $(x,y) \in \lambda$ and consider the sequence $t(1,1) < t(2,1)< \cdots <t(x,1) <t(x,2) < \cdots <t(x,r+1) < t(x+1,r+1) < \cdots <t(g-d+r,r+1)$.
It implies $\im (t)$ has at least $g-d+2r$ different values and this implies $\dim (\mathbb{T}(t)) \leq d-2r$.
Also in case $\dim (\mathbb{T}(t))=d-2r$ it implies the sequence is equal to $i_1 < i_2 < \cdots < i_{g-d+2r}$.
This implies $t(x,y)=i_{x+y-1}$ for $(x,y) \in \lambda$.
As in the proof of Proposition \ref{prop1} this implies $m_{i_2}= \cdots = m_{i_{g-d+2r-1}}=2$.
\end{proof}

In order to prove Theorem A we also need the following lemma. Because of this lemma in the proof of Theorem A, we are going to be able to find a lot of $\underline{m}$-displacement tableaux as in Lemma \ref{lemma2} .

\begin{lemma}\label{lemma3}
Let $\Gamma$ be a chain of cycles of genus $g$ and let $D$ be an effective divisor on $\Gamma$ such that its restriction to some $C_i$ with $2 \leq i \leq g-1$ is equal to $\langle\xi \rangle_i$ with $\xi \notin \mathbb{Z}+\mathbb{Z}\cdot[\ell_i/\ell(v_i,w_i)]$.
In case $\rank{D}=r\geq 1$ then each $\underline{m}$-displacement tableau $t$ on $\lambda=[(g-d+r) \times (r+1)]$ such that the representing divisor of $D$ belongs to $\mathbb{T}(t)$ satisfies $i \notin \im(t)$.
\end{lemma}

\begin{proof}
Let $\sum_{j=1}^g \langle\xi _j\rangle_j +(d-g)w_g$ be the representing divisor of $D$.
If there is an $\underline{m}$-displacement tableau $t$ representing this divisor such that $i \in \im (t)$ then $\xi _i \in \mathbb{Z}$.
There exists a rational function $f$ on $\Gamma$ such that $\divi (f)=D-\sum_{j=1}^g \langle\xi _j\rangle_j - (d-g)w_g$.
Restricting $f$ to $C_i$ we obtain a rational function $f_i$ on $C_i$ such that $\divi (f_i)=a.v_i+b.w_i + \langle\xi \rangle_i - \langle\xi _i\rangle_i$ for some $a,b \in \mathbb{Z}$.

Since $\xi_i \in \mathbb{Z}$ this implies $\xi \in \mathbb{Z}+\mathbb{Z}\cdot[\ell_i/\ell(v_i,w_i)]$.
To explain this claim assume the situation is as shown in Figure \ref{Figuur 5}

\begin{figure}[h]
\begin{center}
\includegraphics[height=2 cm]{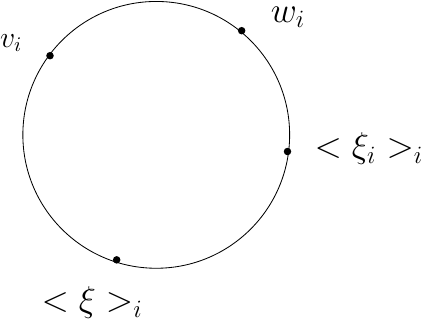}
\caption{Illustration of the proof of Lemma \ref{lemma3} }\label{Figuur 5}
\end{center}
\end{figure}

Up to changing $\xi _i$ and $\xi$ by an integer multiple of $\ell_i/\ell(v_i,w_i)$ we can assume that $0< \xi _i < \xi <(\ell_i-\ell(v_i,w_i))/\ell(v_i,w_i)$.
By definition $f_i$ has integer slopes $c_1$ on the positive oriented arc from $v_i$ to $w_i$, $c_2$ from $w_i$ to $\langle\xi _i\rangle_i$, $c_3$ from $\langle\xi _i\rangle_i$ to $\langle\xi \rangle_i$ and $c_4$ from $\langle\xi \rangle_i$ to $v_i$ (with e.g. $c_4 - c_3=1$).
One has $f(v_i)=f(v_i)+(c_1+c_2\cdot\xi_i+c_3\cdot(\xi - \xi_i)+c_4\cdot(\ell_i/\ell(v_i,w_i)-(\xi +1))\cdot\ell(v_i,w_i)$.
This implies $\xi = c_1+\xi_i+c_4\cdot(\ell_i/\ell(v_i,w_i)-1)$.
Similar computations can be made in case the distribution of the points is different from Figuur \ref{Figuur 5}.

So we obtain a contradiction to the assumptions, hence $i \notin \im(t)$.
\end{proof}

Using the notation of the previous lemma, for the existence of $\langle\xi \rangle_i$ implying $i \notin \im(t)$ one can also argue as follows.
The divisors $av_i+bw_i-\langle\xi_i\rangle_i$ with $a,b,\xi_i \in \mathbb{Z}$ and $a+b=2$ are a countable set of divisors of degree 1 on $C_i$.
They are equivalent to $\langle\xi(a,b)\rangle_i$ for some $\xi (a,b) \in \mathbb{R}$.
This gives a countable subset of $\mathbb{R}$ and we can assume $\xi$ does not belong to it.

A combination of Lemma \ref{lemma2} and \ref{lemma3} gives rise to the following main ingredient of the proofs in this section.

\begin{lemma}\label{extra1}
Let $\Gamma$ be a chain of cycles of genus $g$, let $r,d$ be integers with $1 \leq r \leq g-2$ and $2r \leq d \leq g-2+r$ and assume $w^r_d(\Gamma)=d-2r$.
Fix $J \subset \{ 1, \cdots , g\}$ with $\vert J \vert=d-r$.
Then there exists $J' \subset J$ with $\vert J' \vert=d-2r$ such that all elements $i$ in the complement of $J'$ in $\{ 1, \cdots , g\}$ not being the smallest or largest one satisfy $m_i=2$.
\end{lemma}

\begin{proof}
For $j \in J$ choose $\langle\xi _j\rangle_j \in C_j$ with $\xi _j \notin \mathbb{Z}+\mathbb{Z}\cdot[\ell_j/\ell(v_j,w_j)]$.
Since $w^r_d(\Gamma)=d-2r$ there exists an effective divisor $D$ of degree $d$ containing the sum of the points $\langle\xi _j\rangle_j$ for $j\in J$ and such that $rk (D)\geq r$. 
Since $D$ contains only $r$ points besides the points $\langle\xi _j\rangle_j$ with $j \in J$, we can find $J' \subset J$ consisting of $d-2r$ elements such that $D$ restricted to $C_j$ for $j \in J'$ is equal to $\langle\xi _j\rangle_j$.

The representing divisor of $D$ belongs to $\mathbb{T}(t)$ with $t$ an $\underline{m}$-displacement tableau on $[(g-d+r) \times (r+1)]$.
From Lemma \ref{lemma3} we know $j \notin \im (t)$ for $j \in J'$.
This implies $\dim(\mathbb{T}(t))\geq d-2r$.
Then from Lemma \ref{lemma2} we know that $\im (t)$ is the complement of $J'$ in $\{ 1, \cdots , g\}$ and each $i \in \im (t)$ not being the smallest or largest one satisfies $m_i=2$.
\end{proof}

Now we can prove Theorem A from the introduction.
\begin{theorem}\label{extraTh1}
Let $\Gamma$ be a chain of cycles of genus $g$.
Let $r$ be an integer with $1 \leq r \leq g-2$ and let $d$ be an integer with $2r \leq d \leq g-3+r$.
Then $w^r_d(\Gamma)=d-2r$ implies $\Gamma$ is hyperelliptic.
\end{theorem}

\begin{proof}
We need to prove that for each $2 \leq i \leq g-1$ one has $m_i=2$.
Take a  subset $J$ from $\{ 1, \cdots, g \}$ not containing any element from $\{i-1, i, i+1 \}$ consisting of $d-r$ numbers.
This is possible since $d-r+3 \leq g$.

Lemma \ref{extra1} implies the existence of $J' \subset J$ with $\vert J' \vert=d-2r$ such that each $j$ in the complement of $J'$ in $\{ 1, \cdots ,g \}$ different from the smallest or largest one satisfies $m_j=2$.
Since no element of $\{ i-1,i,i+1 \}$ belongs to $J'$ we find $m_i =2$.
\end{proof}

Note that in the previous proof it is not sufficient to have $i \notin J'$ because we need $i \neq i_1$ and $i \neq i_{g-d+2r}$ using the notations of Lemma \ref{lemma2}.

In Lemma \ref{lemma6} we prove the existence of non-hyperelliptic chains of cycles $\Gamma$ of genus $g$ with $\dim (W^1_d(\Gamma))=d-2$ for some $d \leq g-2$.
This proves that, also when restricting to chains of cycles, Theorem \ref{extraTh1} does not hold when $w^r_d (\Gamma)$ is replaced by $\dim (W^r_d (\Gamma))$.

As mentioned in the introduction, comparing Theorem A with the statement of H. Martens' Theorem for curves, we are missing the cases $(d,r)=(g+r-2,r)$ for $1 \leq r \leq g-2$.
In those cases, if $w^r_{g-2+r}(\Gamma)=g-r-2$, it is possible, using the arguments of the proof of Theorem A, to get some conditions on $\Gamma$.

\begin{lemma}\label{lemma4}
Assume $\Gamma$ is a chain of cycles of genus $g$ and $w^r_{g+r-2}(\Gamma)=g-2-r$ for some integer $1 \leq r \leq g-2$.
Then $m_i=2$ for $2 \leq i \leq 1+r$ and for $g-r \leq i \leq g-1$.
\end{lemma}

\begin{proof}
Assume $2 \leq i \leq r+1$ and take $J=\{ 2, \cdots , g \}\setminus \{ i \}$.
From Lemma \ref{extra1} we find there exist $j_1 < \cdots < j_r$ belonging to $J$ such that each $j \in \{ 1, i, j_1, \cdots , j_r \}$ not being the smallest or largest one satisfies $m_j=2$.
Since $i \leq r+1$ it follows $j_r > i$.
It follows that $m_i=2$.
The proof for $g-r\leq i \leq g-1$ is similar.
\end{proof}

\begin{corollary}
If $\Gamma$ is a chain of cycles of genus $r+2 \leq g \leq 2r+2$ with $r \geq 1$ and $w^r_{g+r-2}=g-r-2$ then $\Gamma$ is hyperelliptic.
\end{corollary}

\begin{lemma}\label{lemma5}
Let $\Gamma$ be a chain of cycles of genus $g\geq 2r+3$ for some $r\geq 1$ such that $w^r_{g+r-2}(\Gamma)=g-r-2$.
Assume $m_i \neq 2$ for some $r+2 \leq i \leq g-r-1$ then $m_{i-l}=m_{i+l}=2$ for $1 \leq l \leq r$.
\end{lemma}

\begin{proof}
Take some $1 \leq l \leq r$ and take $J=\{1, \cdots ,g\} \setminus \{ i-l,i \}$.
Because of Lemma \ref{extra1} there exist $j_1 < \cdots < j_r$ belonging to $J$ such that each $j \in \{ i-l, i, j_1, \cdots ,j_r \}$ not being the smallest or largest one satisfies $m_j=2$.
In case $i-l < j_1$ then we have $i<j_r$ and we find $m_i=2$, which is a contradiction.
So we find $j_1 <i-l$ and therefore $m_{i-l}=2$.
The proof for $m_{i+l}=2$ is similar.
\end{proof}
The remaining non-hyperelliptic chains of cycles of genus $g \geq 2r+3$ that could a priori satisfy $w^r_{g+r-2}(\Gamma)=g-2-r$ are of a very special type. 
We introduce the following definition, describing those chains of cycles.

\begin{definition}\label{def5}
A chain of cycles of genus $g \geq 2r+3$ is called Martens-special for rank $r$ and of type $k$ for some integer $k \geq 1$ if there exist $r+1 < j_1 < \cdots < j_k <g-r$ with $j_{i+1}-j_i \geq r+1$ for each $1 \leq i \leq k-1$ such that the following conditions hold on the torsion profile of $\Gamma$.
For $1 \leq i \leq k$ we have $m_{j_i} \neq 2$.
For $2 \leq i \leq g-1$ with $i \notin \{j_1, \cdots ,j_k \}$ we have $m_i=2$.
\end{definition}

Sometimes we omit mentioning the type of a Martens-special chain of cycles for some rank $r$.
As already mentioned in the introduction, in case $r \geq 2$, then a Martens-special chain of cycles for rank $r$ is also a Martens-special chain of cycles for rank $r-1$.
Therefore, in case $r=1$, we also talk about Martens-special chains of cycles.

From Lemmas \ref{lemma4} and \ref{lemma5} we obtain the ``only if'' part of Theorem B from the introduction.
\begin{theorem}
Let $\Gamma$ be a chain of cycles of genus $g$ and let $r$ be an integer with $1 \leq r \leq g-2$.
If $w^r_{g+r-2}(\Gamma)=g-r-2$ then either $\Gamma$ is hyperelliptic or $\Gamma$ is Martens-special for rank $r$.
\end{theorem}

In the next section we are going to prove that those Martens-special chains of cycles for rank $r$ do satisfy $w^r_{g-2+r}(\Gamma)=g-2-r$, hence completing the proof of Theorem B from the introduction.
So, in the case of chains of cycles, we obtain counterexamples to the statement of H. Martens' theorem as for the case of curves.
In particular, in case $r=1$, we obtain counterexamples to the conjecture in \cite{ref10}. 
Within the set of chains of cycles we then classified all such counterexamples.

\section{Non-hyperelliptic chains of cycles satisfying $w^r_{g-2+r}=g-2-r$ for some $1 \leq r \leq g-2$}\label{section4}

In the proof of Theorem B we are going to make use of the tropical Abel Theorem obtained in e.g. \cite{refMZ}, Theorem 6.1.
We refer to \cite{refMZ} for the general definitions and now we give a description in the case of a chain of cycles $\Gamma$.

Associated to $\Gamma$ we have a space of 1-forms $\Omega (\Gamma)$ and for $\omega \in \Omega (\Gamma)$ and $\gamma$ a path in $\Gamma$ we have an integral $\int _{\gamma}\omega \in \mathbb{R}$.
In case $\Gamma$ is a single cycle $C$ and we choose an orientation on $C$ then there is a unique $\omega _C\in \Omega(C)$ such that $\int _{\gamma} \omega _C$ is the arclength of $\gamma$ (with a sign according to the orientation).
For $1 \leq i \leq g$ let $\omega_i$ be the 1-form on $\Gamma$ obtained by extending $\omega _{C_i} \in \Omega (C_i)$ by 0 outside $C_i$.
We have $\Omega (\Gamma) = \oplus_{i=1}^g \mathbb{R}\omega_i$.

For $1 \leq i \leq g$ we have $\int _{\Gamma} \omega _i=\ell_i$.
In writing this integral, we consider $\Gamma$ as a path, starting at $w_1$ and following each $C_i$ with the positive orientation.
Let $\{ \omega _i^* \}$ be the dual base of $\{ \omega _i \}$ in the dual space $\Omega (\Gamma)^*$.
The torus $\Pi _{i=1}^g (\mathbb{R} \omega _i^* / \mathbb{Z} \ell_i\omega_i^*)$ is the Jacobian $J(\Gamma)$ of $\Gamma$ (it can be identified with $\Pi_{i=1}^g C_i$).
Then for $p_1, p_2 \in \Gamma$ and using a path $\gamma$ from $p_1$ to $p_2$, the linear map $l_{\gamma} : \Omega (\Gamma) \rightarrow \mathbb{R} : \omega \rightarrow \int _{\gamma} \omega$ defines an element of $J(C)$ only depending on $p_1$ and $p_2$ and not on the choice of $\gamma$.
We write $\int _{p_1}^{p_2} \omega _j \in \mathbb{R} / \mathbb{Z} \ell_j$ to denote the projection of $l_{\gamma}$ on $\mathbb{R} \omega _j^* / \mathbb{Z}\ell_j \omega_j^* \cong \mathbb{R} / \mathbb{Z} \ell_j$.
For a divisor $E = \sum _{i=1}^d P_i$ of degree $d$ on $\Gamma$ we write $\int _{d\cdot w_g}^E \omega_j=\sum_{i=1}^d \int _{w_g}^{P_i} \omega _j$.
Let $S^d (\Gamma)$ be the space of effective divisors of degree $d$ on $\Gamma$ (a symmetric product of $\Gamma$).
The tropical Abel-map is $I(d) : S^d(\Gamma) \rightarrow J(\Gamma) : E \rightarrow (\int _{d\cdot w_g}^E \omega_j)_{j=1}^g$.

\begin{theorem}[Tropical Abel Theorem]\label{TheoremTAT}
If $E_1, E_2 \in S^d(\Gamma)$ then $E_1$ is equivalent to $E_2$ if and only if $I(d)(E_1)=I(d)(E_2)$.
\end{theorem}

We now give the main lemma based on the tropical Abel Theorem that will be used in the proof of Theorem B.
This lemma is already noticed in \cite{JR}, 2.2.

\begin{lemma}\label{lemmaTAT}
Let $E$ be an effective divisor of degree $d$ on $\Gamma$ and let $\sum _{i=1}^g \langle\xi _i\rangle_i$ be the representing divisor.
For $1 \leq j \leq g$ the point $\langle\xi _j\rangle_j$ on $C_j$ is determined by the equation $\int _{w_j}^{\langle \xi _j\rangle_j} \omega _j \equiv \int _{d\cdot w_g}^E \omega _j +(j-1)\ell(v_j,w_j)$ in $\mathbb{R} / \mathbb{Z} \ell_j$.
\end{lemma}

\begin{proof}
By definition of the representing divisor the difference of $E$ and $\sum _{i=1}^g \langle\xi _i \rangle_i$ is a multiple of $w_g$.
From Theorem \ref{TheoremTAT} it follows that for each $1 \leq j \leq g$ we have $\int _{d\cdot w_g}^E \omega _j \equiv \sum _{i=1}^g \int _{w_g}^{\langle\xi _i\rangle_i} \omega _j$.
For $i>g$ we have $\int _{w_g}^{\langle\xi _i\rangle_i} \omega _j=0$ because there is a path of $w_g$ to $\langle\xi _i\rangle_i$ not meeting $C_j$.
For $i<j$ there is a path from $w_g$ to $w_i$ going from $w_j$ to $v_j$ negatively oriented on $C_j$, hence $\int _{v_g}^{\langle\xi _i\rangle_i} \omega _j \equiv -\ell(v_j,w_j)$.
Of course $\int _{w_g}^{\langle\xi _j\rangle_j} \omega _j \equiv \int _{w_j}^{\langle\xi _j\rangle_j} \omega _j$.
This implies the lemma.
\end{proof}

Note that Lemma  \ref{lemmaTAT} gives a very short proof of Lemma \ref{lemma1} using the tropical Abel Theorem.

The next theorem finishes the proof of Theorem B from the introduction.
\begin{theorem}\label{theoremRB}
A Martens-special chain of cycles $\Gamma$ for rank $r$ satisfies $w_{g-2+r}^r(\Gamma)=g-2-r$.
\end{theorem}

\begin{proof}
As mentioned in the introduction, using the arguments from the appendix of \cite{ref10}, we have that $w^r_{g-2+r}(\Gamma) \leq g-2-r$.
Therefore from Definition \ref{def3} it is enough to prove that for each $P_1, \cdots , P_{g-2}$ on $\Gamma$ there exists an effective divisor $E$ of degree $g-2+r$ and rank at least $r$ containing $P_1 + \cdots + P_{g-2}$.
To prove this it is allowed to replace $P_1 + \cdots + P_{g-2}$ by any equivalent effective divisor (see the argument before Lemma \ref{lemma1}).
Using equivalence of divisors on chains of cycles we can therefore assume that there exist $1 \leq i_1 < i_2 < \cdots < i_{g-2} \leq g$ such that $P_j \in C_{i_j}$ for all $1 \leq j \leq g-2$.

Let $a<b$ with $\{ 1, \cdots , g \} \setminus \{i_1 , \cdots , i_{g-2} \}=\{ a,b \}$.
We distinguish between 4 cases depending on whether $a,b$ belong to $\{ j_1, \cdots , j_k \}$ or not (remember the notation introduced in Definition \ref{def5}).
In each case we give a certain divisor $E$ of degree $g-2+r$ containing $P_1 + \cdots + P_{g-2}$ and we prove $\rank (E) \geq r$.

\begin{enumerate}
\item Assume there exists $1 \leq i <i' \leq k$ such that $a=j_i$ and $b=j_{i'}$.

For $1 \leq e \leq r$, since $j_{i+1}-j_i=j_{i+1}-a \geq r+1$, we have $m_{a+e}=2$.
We take $P'_{a+e} \in C_{a+e}$ such that $P_{a+e}+P'_{a+e}$ is equivalent to $2v_{a+e}$ (hence also to $2w_{a+e}$) and we consider $E = \sum_{i=1}^{g-2} P_i + \sum_{e=1}^r P'_{a+e}$.
Let $D=\sum _{i=1}^g \langle\xi _i\rangle_i + (r-2) w_g$ be the representing divisor of $E$.

Since $E$ contains no point on $C_b$ and exactly $b-2+r$ points on the cycles $C_i$ with $i<b$ we have $\int_{(g+r-2)w_g}^E \omega_b \equiv -(b+r-2)\ell(v_b,w_b)$.
Because of Lemma \ref{lemmaTAT} this implies $\int _{w_b}^{\langle\xi _b\rangle_b} \omega _b \equiv -(r-1)\ell(v_g,w_g)$.
This implies $<\xi _b>_b = <1-r>_b$.

For $1 \leq e \leq r$ we have $E$ contains exactly $P_{a+e}+P'_{a+e}$ on $C_{a+e}$.
Since $P_{a+e}+P'_{a+e}$ is equivalent to $2w_{a+e}$ on $C_{a+e}$, we have $\int _{2w_{a+e}}^{P_{a+e}+P'_{a+e}} \omega _{a+e}= \int _{2w_{a+e}}^{2w_{a+e}} \omega _{a+e} \equiv 0$.
Also there are exactly $a+2e-3$ points of $E$ on the cycles $C_i$ with $i<a+e$.
This implies $\int _{(g+r-2)w_g}^E \omega _{a+e} = -(a+2e-3)\ell(v_{a+e}, w_{a+e})$ and therefore from Lemma \ref{lemmaTAT} it follow $\int _{w_{a+e}}^{\langle\xi _{a+e} \rangle_{a+e}} \omega_{a+e} \equiv -(e-2)\ell(v_{a+e},w_{a+e})$.
Since $m_{a+e}=2$ this means $\langle\xi _{a+e} \rangle_{a+e}=v_{a+e}$ if $e$ is odd and $\langle\xi _{a+e}\rangle_{a+e}=w_{a+e}$ if $e$ is even.

There is no point of $E$ on $C_a=C_{j_i}$ and exactly $a-1$ points of $E$ on the cycles $C_i$ with $i<a$.
This implies $\int _{(g+r-2)w_g}^E \omega _a \equiv -(a-1)\ell(v_a,w_a)$.
Because of Lemma \ref{lemmaTAT} we need $\int _{w_a}^{\langle\xi_a\rangle_a} \omega _a \equiv 0$, hence $\langle\xi _a\rangle_a=w_a$.

On $[2 \times (r+1)]$ consider the tableau $t$ defined by $t(i,j)=a+i+j-2$ if $(i,j)\neq (2,r+1)$ and $t(2,r+1)=b$.
Since $m_{a+i}=2$ for $1 \leq i \leq r$ this is an $\underline{m}$-displacement tableau.
Moreover for $i=t(x,y)$ we found $\langle\xi _i\rangle_i=\langle x-y\rangle_i$, hence $D \in \mathbb{T}(t)$.
Because of Theorem \ref{theorem1} this implies $\rank (D)= \rank (E) \geq r$.

Note that if $i'=i+1$ we really need $j_{i+1}-j_i \geq r+1$ for this part of the proof.

\item Assume $m_a=2$ or $a=1$ and $b=j_i$ for some $1 \leq i \leq k$.

In the proof we are going to distinguish two subcases: $a \geq r+1$ and $a < r+1$.

In both cases, for $2 \leq i \leq r+1$ we have $m_i=2$.
In case $a \geq r+1$ then for $1 \leq i \leq r$ we take $P'_i \in C_i$ such that $P_i + P'_i$ is equivalent to $2w_i$ (hence also to $2v_i$ in case $i>1$) and in case $a<r+1$ then for $1 \leq i \leq r+1$ with $i \neq a$ we take $P'_i \in C_i$ satisfying the same condition.
We define $E$ as being the sum of the points $P_i$ and $P'_i$ and we write $D= \sum_{i=1}^g \langle\xi _i\rangle_i + (r-2)w_g$ for the representing divisor of $E$.

In both subcases, for $i <r+1$ with $i<a$, the divisor $E$ contains exactly $P_i+P'_i$ on $C_i$ and exactly $2(i-1)$ points on the cycles $C_j$ with $j<i$.
This implies $\int _{(g+r-2)w_g}^E \omega _i \equiv -(2(i-1))\ell(v_i,w_i)$ and therefore because of Lemma \ref{lemmaTAT} $\int _{w_i}^{\langle\xi _i\rangle_i} \omega _i \equiv -(i-1)\ell(v_i;w_i)$.
It follows $\langle\xi _i\rangle_i=v_i$ in case $i$ is even and $\langle\xi _i\rangle_i=w_i$ in case $i$ is odd.

In case $a \geq r+1$ there is no point of $E$ on $C_a$ and there are exactly $a+r-1$ points of $E$ on the cycles $C_j$ with $j<a$.
This implies $\int _{(g+r-2)w_g}^E \omega _a \equiv -(a-1+r)\ell(v_a,w_a)$ and because of Lemma \ref{lemmaTAT} we obtain $\int _{w_a}^{\langle\xi _a\rangle_a} \omega _a \equiv -r\cdot \ell(v_a,w_a)$.
Since $m_a=2$ this implies $\langle\xi_a\rangle_a =w_a$ if $r$ is even and $\langle\xi _a\rangle_a= v_a$ if $r$ is odd.

Using exactly the same argument as in the previous case, we also have $\langle\xi _b\rangle_b=\langle1-r\rangle_b$.

On $[2 \times (r+1)]$ we define $t(2,r+1)=b$, $t(1,r+1)=t(2,r)=a$ and $t(i,j)=i+j-1$ otherwise.
Since $m_2 = \cdots = m_r = m_a=2$ this is an $\underline{m}$-displacement tableau.
Since $D \in \mathbb{T}(t)$ we obtain $\rank (D)=\rank (E) \geq r$.

Now assume $a<r+1$.
We leave it to the reader to check (in the same way as before) that also $\langle\xi _i\rangle_i=v_i$ if $i$ is even and $\langle\xi _i\rangle_i=w_i$ if $i$ is odd for $a \leq i \leq r+1$ and we still have $\langle\xi _b\rangle_b=\langle1-r\rangle_b$.

Now on $[2 \times (r+1)]$ we define $t(2,r+1)=b$ and $t(i,j)=i+j-1$ otherwise.
Again this is an $\underline{m}$-displacement tableau and $D \in \mathbb{T}(t)$, showing $\rank (D)=\rank (E) \geq r$.

Note that in this second subcase, if $i=1$ we really need $j_1 \geq r+2$ for this part of the proof.

\item Assume there exists $1 \leq i \leq k$ such that $a=j_i$ and $m_b=2$ or $b=g$.

This possibility is the same as the previous one replacing $C_i$ by $C_{g-i}$, $v_i$ by $w_{g-i}$ and $w_i$ by $w_{g-i}$ for $1 \leq i \leq g$.
In particular we use the conditions $m_i=2$ for $g-r \leq i \leq g-1$ and we really need this condition if $i=k$.

\item Assume $a=1$ or $m_a=2$ and $b=g$ or $m_b=2$.

Let $i_1 < i_2 < \cdots < i_r$ be the $r$ smallest integers from $\{ 1, 2, \cdots , g\} \setminus \{ j_1, \cdots , j_k, a, b \}$ and for $i \in \{ i_1, \cdots ,i_r \}$ take $P'_i \in C_i$ such that $P_i+P'_i$ is equivalent to $2 w_i$.
Take $E= \sum _{i=1}^{g-2}P_i + \sum _{j=1}^r P'_{i_j}$ and let $D = \sum _{i=1}^g \langle\xi _i\rangle_i +(r-2)w_g$ be the representing divisor.
We have to distinguish between three subcases: $a \geq r+1$; $a \leq r$ and $b>r+1$; $b \leq r+1$ (remember $a<b$).

In all subcases, as in case (2), if $i \leq r$ with $i \leq a$ we find: $\langle\xi _i\rangle_i=w_i$ if $i$ is odd and $\langle\xi _i\rangle_i=v_i$ is $i$ is even.

Assume $a \geq r+1$.
Then as in the case (2) we find $\langle\xi _a\rangle_a=w_a$ if $r$ is even and $\langle\xi _a\rangle_a=v_a$ if $r$ is odd.
Also we still obtain $\langle\xi _b\rangle_b=\langle1-r\rangle_b$.
Then we use exactly the same $\underline{m}$-tableau on $[2 \times (r+1)]$ as in the first subcase of case (2) to conclude $\rank (D)=\rank (E)\geq r$.

Next assume $a \leq r$ but $b>r+1$.
As in the second subcase of case (2) for $1 \leq i \leq r+1$ we find: $\langle\xi _i\rangle_i=w_i$ if $i$ is odd and $\langle\xi _i\rangle_i=v_i$ if $i$ is even.
We also obtain $\langle\xi _b\rangle_b = \langle1-r\rangle_b$.
So using exactly the same $\underline{m}$-displacement tableau as in the second subcase of case (2) we obtain $\rank (D)=\rank (E) \geq r$.

So finally assume $b \leq r+1$.
In case $j_1 >r+2$ we have $\{ i_1, \cdots , i_r,a,b\} =\{ 1,2, \cdots , r+2 \}$ and we find $\langle\xi _i\rangle_i=v_i$ if $i$ is even and $\langle\xi _i\rangle_i=w_i$ if $i$ is odd for $1 \leq i \leq r+2$.
On $[2 \times (r+1)]$ we use $t(i,j)=i+j-1$.
It is an $\underline{m}$-displacement tableau and $D \in \mathbb{T}(t)$, this proves $\rank (D)=\rank (E) \geq r$.

In case $j_1=r+2$ then $\{ i_1, \cdots , i_r, a, b \}=\{ 1, \cdots, r+1,r+3 \}$.
The divisor $E$ contains exactly $P_{r+3}+P'_{r+3}$ on $C_{r+3}$ (and remember $m_{r+3}=2$) and it has exactly $2r-1$ points on the cycles $C_j$ with $j<r+3$.
This implies $\int _{(g+r-2)w_g}^E \omega _{r+3} \equiv -(2r-1)\ell(v_{r+3},w_{r+3})$ and therefore because of Lemma \ref{lemmaTAT} $\int _{w_{r+3}}^{\langle\xi _{r+3}\rangle_{r+3}} \omega_{r+3} \equiv (3-r)\ell(v_{r+3},w_{r+3})$.
Since $m_{r+3}=2$ we find $\langle\xi _{r+3}\rangle_{r+3} = \langle1-r\rangle_{r+3}$.

For $1 \leq i \leq r+1$ we still find $\langle\xi _i\rangle_i=w_i$ if $i$ is odd and $\langle\xi _i\rangle_i=v_i$ is $i$ is even.

Now on $[2 \times (r+1)]$ we take $t(2,r+1)=r+3$ and $t(i,j)=i+j-1$ otherwise.
Again this is an $\underline{m}$-displacement tableau and we conclude once more that $\rank (D)=\rank (E) \geq r$.

\end{enumerate}

This finishes the proof of Theorem B.
\end{proof}

Combining Theorems A and B we now give the proof of Theorem C.
More concretely we obtain the following theorem.

\begin{theorem}
Let $\Gamma$ be a Martens-special chain of cycles for rank $r$ for some $r \geq 2$.
Then $w_{g-2+r}^r (\Gamma) \neq w_{2g-2-(g-2+r)}^{r-(g-2+r)+g-1} (\Gamma)$.
\end{theorem}

\begin{proof}
From Theorem B we obtain $w^r_{g-2+r}(\Gamma)=g-2-r$.
But, since $\Gamma$ is not hyperelliptic, from Theorem A we obtain $w^{r-(g-2+r)+g-1}_{2g-2-(g-2+r)}(\Gamma)=w^1_{g-r}(\Gamma)<g-r-2$.
This implies the inequality of the theorem.
\end{proof}

\begin{lemma}\label{lemma6}
Let $\Gamma$ be a Martens-special chain of cycles of type $k$.
The graph $\Gamma$ has gonality at most $k+2$ and $\dim (W^1_{k+2}(\Gamma))=k$.
\end{lemma}

\begin{proof}
If $\dim (W^1_{k+2}(\Gamma)) >k$ then there exists an $\underline{m}$-displacement tableau $t$ on $[(g-k-1) \times 2]$ such that the image of $t$ contains at most $g-k-1$ integers.
This is impossible, therefore $\dim (W^1_{k+2}(\Gamma))\leq k$.

Let $\{ i_1, i_2, \cdots , i_{g-k} \}=\{1, \cdots , g\} \setminus \{j_1, \cdots ,j_k \}$ with $i_1 < \cdots < i_{g-k}$.
Take the tableau $t$ on $[(g-k-1) \times (2)]$ with $t(n,1)=i_n$ and $t(n,2)=i_{n+1}$.
Since $m_{i_n}=2$ for $2 \leq n \leq g-k-1$ we obtain that $t$ is an $\underline{m}$-displacement tableau for $\Gamma$.
This proves $\mathbb{T}(t) \subset W^1_{k+2}(\Gamma)$ with $\dim (\mathbb{T}(t))=k$.
\end{proof}

Remember that for a "general" chain of cycles of genus $g$ we have that the gonality is equal to $[(g+3)/2]$ (see \cite{ref9}), so obtaining a smaller gonality for a chain of cycles implies it is of a special type.

\begin{proposition}\label{prop2}[Weak H. Martens' Theorem for $w^r_{g-2+r}$ for chains of cycles]
If $\Gamma$ is a chain of cycles of genus $g$ with $w^r_{g-2+r}(\Gamma)=g-2-r$, then the gonality of $\Gamma$ is at most $(g+r)/(r+1)$ (hence less than $[(g+3)/2])$.
\end{proposition}

\begin{proof}
From Definition \ref{def5} we find that for a Martens-special chain of cycles for rank $r$ of genus $g$ and type $k$ we have $g\geq (r+1)(k+1)+1$.
It follows from Lemma \ref{lemma6} that the gonality is at most $k+2\leq (g+r)/(r+1)$.
\end{proof}

In the forthcoming paper \cite{refCO} we are going to prove that in``general" a Martens-special chain of cycles of type $k$ has gonality equal to $k+2$.
In particular this implies the Weak H. Martens' Theorem for $w^r_{g-2+r}$ for chains of cycles is sharp.
This occurs if $m_j=0$ in case $m_j \neq 2$.
Also despite $\dim (W^1_{k+2}(\Gamma))=k$ we will prove $w^1_{k+2}(\Gamma)=0$ for such general Martens-special chain of cycles of type $k$.
\\

\textbf {Acknowledgement.}  I am grateful to the referee for his suggestions to improve the paper.

\end{document}